\documentclass[preprint, 12pt, 3p]{elsarticle_arxiv}

\usepackage{graphicx}
\usepackage{amssymb}
\usepackage{amsmath}
\usepackage{bm}
\usepackage{a4wide}
\usepackage{appendix}
\usepackage{accents} 

\usepackage{color}
\usepackage{esint}

\newcommand{\red}[1]{{{#1}}}
\newcommand{\blue}[1]{{{#1}}}
\definecolor{lblue}{rgb}{0.0,0.4,1.0}

\definecolor{dgreen}{rgb}{0.0,0.4,0.11}
\newcommand{\dgreen}[1]{{{#1}}}

\newcommand{\magenta}[1]{{{#1}}}

\newcommand{\RELEASE}[1]{}
\newcommand\Appx[1]{\ref{#1}}

\def\Dlt{\Delta}

\newcommand\zerobm{\bf{0}}

\newcommand\dlt{\delta}          
          
\newcommand\vphi{\varphi}

\newcommand\plper{{\ul{\#}}} 

\newcommand\Lpar{[\kern-0.17em{[}}
\newcommand\Rpar{]\kern-0.18em{]}}

\def\FT#1{{\underaccent{{\sim}}{{#1}}}} 

\def\gradpl{\ol{\nabla}}
\def\gradplx{\ol{\nabla}_x}

\def\gradply{\ol{\nabla}_y}

\def\ext{{\rm{ext}}}

\def\vis{{\rm{vis}}}

\newcommand\cwto{\rightharpoonup}

\def\wrt{{w.r.t.{~}}}
\def\rhs{{r.h.s.{~}}}
\def\lhs{{l.h.s.{~}}}
\def\ie{{\it i.e.~}}
\def\eg{{\it e.g.~}}

\newcommand{\eq}[1]{(\ref{#1})}

\def\bmi#1{\textbf{\textit{#1}}}

\newcommand{\ul}[1]{\underline{#1}}
\newcommand{\ol}[1]{\overline{#1}}
\newcommand{\pl}[1]{{#1}'}

\def\pdiff#1#2{\frac{\partial {#1}}{\partial {#2}}}
\newcommand{\pdnw}[1]{\frac{D_{w,\omega}{#1}}{D\,n}} 

\def\taubf{{\mbox{\boldmath$\tau$\unboldmath}}}
\def\nubf{{\mbox{\boldmath$\nu$\unboldmath}}}

\def\R{\hbox{\rm I\kern-0.2em R}}
\def\Z{\hbox{\rm Z\kern-0.3em Z}}
\def\Eop{{{\rm I} \kern-0.2em{\rm E}}}%
\def\Dop{{{\rm I} \kern-0.2em{\rm D}}}%
\def\Cop{{{\rm C} \kern-0.6em{\rm C}}}

\def\veps{\varepsilon}

\def\vkappa{\varkappa}

\def\imu{{\rm{i}}} 

\def\Om{\Omega}
\def\om{\omega}
\def\pd{\partial}

\def\RR{{\mathbb{R}}}
\def\ZZ{{\mathbb{Z}}}

\def\ab{{\bmi{a}}}
\def\bb{{\bmi{b}}}

\def\fb{{\bmi{f}}}

\def\nb{{\bmi{n}}}

\def\ub{{\bmi{u}}}

\def\wb{{\bmi{w}}}

\def\Ab{{\bmi{A}}}
\def\Bb{{\bmi{B}}}

\def\Qb{{\bmi{Q}}}

\def\Wb{{\bmi{W}}}

\def\Tcal{\mathcal{T}}

\def\ipYs#1#2{\left ({#1},\,{#2}\right )_{Y^*}}

\def\ipYf#1#2{{\left({#1},\,{#2}\right)_{Y^*}}}

\def\AA{\mathcal{A}_{\wb}}

\def\Tuf#1{{\mathcal{T}}_\veps{\left ({#1}\right )}} 
\newcommand\Tuftxt{\Tcal_\veps\,}

\def\intY{\fint}

\def\Eop{{{\rm I} \kern-0.2em{\rm E}}}%

\newcounter{theorem}
               {\vspace{7pt} \noindent{\bf Theorem
                   \refstepcounter{theorem}\thetheorem. \protect\label{#1}
                   \hspace{-2mm}}\it}%
               {\vspace*{7pt}}

\newcounter{proposition}
               {\vspace{7pt} \noindent{\bf Proposition
                   \refstepcounter{proposition}\theproposition. \protect\label{#1}
                   \hspace{-2mm}}\it}%
               {\vspace*{7pt}}
               
\newcounter{lemma}
               {\vspace{7pt} \noindent{\bf Lemma
                   \refstepcounter{lemma}\thelemma. \protect\label{#1}
                   \hspace{-2mm}}\it}%
               {\vspace*{7pt}}

\newcounter{remark}
\newenvironment{myremark}[1]%
               {\vspace{7pt} \noindent{\bf Remark
                   \refstepcounter{remark}\theremark. \protect\label{#1} \hspace{-2mm}}\rm }%
               {\vspace*{-7pt} \flushright $\triangle$\\ } 

\newenvironment{myremark-noendsign}[1]%
              {\vspace{7pt} \noindent{\bf Remark
                  \refstepcounter{remark}\theremark. \protect\label{#1} \hspace{-2mm}}\rm }%
               {\vspace*{-7pt} \flushright $\square$\\ } 

\newcounter{definition}
               {\vspace{7pt} \noindent{\bf Definition
                   \refstepcounter{definition}\thedefinition. \protect\label{#1} \hspace{-2mm}}\rm }%
               {\vspace*{-7pt} \flushright $\square$\\ } 

\begin{document}

\begin{frontmatter}

\title{Homogenization of the vibro--acoustic transmission on periodically perforated elastic plates interacting with flow}

\author[NTIS]{E.~Rohan\corref{cor1}}
\ead{rohan@kme.zcu.cz}
\cortext[cor1]{Corresponding author}
\author[NTIS]{V.~Luke\v{s}}
\ead{vlukes@kme.zcu.cz}
\address[NTIS]{Department of Mechanics, NTIS New Technologies for Information Society, \\
Faculty of Applied Sciences, University of West Bohemia, \\
Univerzitn\'\i~8, 30100 Pilsen, Czech Republic}

\begin{abstract}
  We consider acoustic waves propagating in an inviscid fluid
  interacting with a rigid periodically perforated plate in the
  presence of permanent flows. The paper presents a model of an
  acoustic interface obtained by the asymptotic homogenization of a
  thin transmission layer in which the plate is embedded.  To account
  for the flow, a decomposition of the fluid pressure and velocity in
  the steady and fluctuating parts is employed. This enables for a
  linearization and an efficient use of the homogenization method
  which leads to a model order reduction effect. The dependence of an
  extended Helmholtz equation on the permanent flow introduces a
  locally periodic velocity field in the perforated plate vicinity, so
  that the coefficients of the homogenized interface depend on the
  flow. The derived model extended by natural coupling conditions
  provides an implicit Dirichlet-to-Neumann operator.  Numerical
  simulations of wave propagation in a waveguide illustrate the flow
  speed influence on the acoustic transmission. Also some geometrical
  aspects are explored.

\end{abstract}


\end{frontmatter}

\section{Introduction}\label{sec-intro}

The problem studied in this paper is motivated by various industrial
applications where designed structures incorporate perforated plates, or panels
which enable for cooling, or ventilation by fluid flow through the apertures,
and simultaneously should reduce the noise transmission. The permanent flow can
be a quite important phenomenon influencing the acoustic field.

We consider fluid acoustics in a waveguide in which a rigid periodically perforated plate is fitted. As the new ingredient of the modelling, the steady fluid flow is respected.
The aim 
is to derive a homogenized model of acoustic waves propagating in a thin layer $\Om_\dlt \subset \RR^3$ occupied by a fluid interacting with the plate. 
The layer is embedded in a waveguide where it separates  two fluid-filled subdomains. 
The derived homogenized model constitutes interface transmission conditions on the ``homogenized acoustic metasurface'' $\Gamma_0 \subset \RR^2 \times \{0\}$ replacing the problem of the fluid acoustics in a complex 3D geometry of the periodically perforated plate. This modelling approach was proposed in \citep{rohan-lukes-waves07}, and further elaborated in \cite{Rohan-Lukes-AMC2019,Rohan-Lukes-AMM2022} for elastic plates. Its usefulness for efficient solving optimum design problems was shown in \cite{noguchi2020topology,Noguchi-Yamada-2021}.
To account for the advection effects related to the permanent flow which is assumed to be independent of the acoustic perturbations, a linearization is employed to establish approximate models of acoustic waves. For further treatment by the homogenization method, the pressure formulation derived in \cite{Rohan-Cimrman-AMC2021} while respecting the advection in an inviscid barotropic fluid is employed.
Up to our knowledge, the acoustic transmission in a nonstationary fluid flowing through a periodically perforated interfaces has not been treated by the homogenization so far in the published literature. Numerical aspects of the acoustic waves in a uniform flow were reported in \cite{Mercier-Maurel-2016} where the applied  Lorentz transformation yields the Helmholtz equation. Another study \cite{BBD-GalbrunEq-2008} employed the Galbrun equation.

Our paper contributes to the research in the two-scale modelling of periodically heterogeneous interfaces. In our previous studies \citep{rohan-lukes-waves07,Rohan-Lukes-AMC2019,Rohan-Lukes-AMM2022} we were concerned with waves in standing fluid. For this, the homogenization strategy has been applied  in situations when thin rigid, or elastic  perforated plate represented by interface $\Gamma_0$ is  characterized by the thickness $\approx \dlt$ which is  proportional to the size of the perforating holes $\veps\approx\dlt$.
 We developed homogenized models of a layer of the thickness $\dlt$, containing a perforated plate. Using an approximation respecting a given finite scale $\veps_0 > 0$, non-local transmission conditions of the acoustic field interacting with the rigid, or elastic plate were obtained as the two-scale homogenization limit $\veps \rightarrow 0$. The same modelling strategy is pursued in this paper.

Analogous problems with thin perforated interfaces have been studied using different approaches in a number of works, \eg   \citep{bb2005,Schweizer-JMPA2020} showing that the first order homogenization $\veps \rightarrow 0$ leads to totally transparent interfaces without any effects of the finite thickness. Besides the acoustic problems with a standing fluid, similar treatment was employed to study the electromagnetic field \cite{DELOURME201228}.
Using higher order approximation involving the correctors at order $o(\veps^1)$,  nontrivial interface conditions capturing acoustic impedance of the thin interfaces have been obtained in 
\cite{Marigo-Maurel-JASA2016b,Pham-Maurel-Marigo-PRSA2021} using an approach based on the so-called inner and outer
asymptotic expansions which enable to treat rather general shapes of the perforations, or other heterogeneities. 

The paper is organized, as follows. The acoustic problem in the waveguide containing a transmission layer is introduced in 
Section~\ref{sec-problem}, where all the geometrical objects are defined and the decomposition of the problem in two-subproblems is established. Then we focus on the acoustic subproblem imposed in the layer; its weak formulation is treated by the asymptotic homogenization in Section~\ref{sec-homog}, where the limit problem is presented, all auxiliary local autonomous corrector problems,  and the macroscopic model of the homogenized acoustic layer are introduced. Expressions for the effective model coefficients involved in the macroscopic model are derived in the \ref{sec-appx1}. In 
Section~\ref{sec:num-simulations}, the homogenized model serving the acoustic transmission condition in the global problem is illustrated  using numerical examples which show the influence of the flow on the acoustic properties of the perforated plate. Concluding remarks and research perspectives follow in Section~\ref{sec-concl}.

 \subsection{Notation}
The spatial position $x$ in the medium is specified through the coordinates $(x_1,x_2,x_3)$ with respect to a Cartesian reference frame $\mathcal{R}(\text{O};\vec{e}_1,\vec{e}_2,\vec{e}_3)$ specified by an orthonormal basis vectors $\vec{e}_k$. The boldface
notation for vectors, $\ab = (a_i)$, and for tensors,  $\bb = (b_{ij})$,
is used. 
 The gradient and divergence operators applied to a vector $\ab$ are
 denoted by $\nabla\ab$ and $\nabla \cdot \ab$, respectively. 
 Alternatively the notation $\nabla_y = (\pd_i^y) = \frac{\pd}{\pd y_i}$ is used in a generic sense.
Throughout the paper, $x$ denotes the global (``macroscopic'') coordinates, while the ``local''
coordinates $y$ describe positions within the representative unit cell
$Y\subset\RR^3$ where $\RR$ is the set of real numbers. By Latin subscripts $i,j,k,l \in\{1,2,3\}$  we refer to vectorial/tensorial components in $\RR^3$. Subscripts $\alpha,\beta \in \{1,2\}$ are reserved for the tangential components with respect to the plate midsurface, \ie coordinates $x_\alpha$ of vector represented by $x' = (x_1,x_2) = (x_\alpha)$ are associated with directions $(\vec{e}_1,\vec{e}_2)$. Moreover, $\gradplx = (\pd_\alpha)$ is the ``in-plane'' gradient. 
By $\nb = (n_i)$ we denote the unit normal vector. The standard notation for Lebesgue $L^k$, and Sobolev $W^{1,2}=H^1,W_0^{1,2}=H_0^1$ functional spaces is adhered.

\section{Problem formulation and decomposition in subproblems}\label{sec-problem}

We find a representation of the acoustic
interaction on a perforated plate under the steady fluid flow. For this, we perform the asymptotic homogenization of the flow and acoustic fields in a transmission layer involving the plate. Desired homogenized acoustic
transmission conditions replacing the layer  are derived using the asymptotic analysis \wrt a scale
parameter $\veps$ which describes the layer thickness and also the size and spacing of holes (general perforations) periodically drilled in the plate structure. In this paper we focus on the homogenization of the acoustic field in the layer. The resulting model is then used to serve the transmission conditions for the interface replacing the layer in the global acoustic problem imposed in a waveguide.

\subsection{Flow and acoustics decomposition}

Following the approach employed in \cite{Rohan-Cimrman-AMC2021} to analyze the acoustic waves propagating in rigid scaffolds, the total fluid fields   
 $\wb$, $p$ and the mass density of the fluid, $\rho$,  are split into the ``stationary flow'' parts $\bar \wb$, $\bar \rho$ and $\bar p$ and the ``acoustic fluctuation'' parts $\tilde \wb$, $\tilde \rho$ and $\tilde p$, so that
\begin{equation}\label{eq-NS2}
   \wb = \bar\wb + \tilde\wb\;,\quad p = \bar p + \tilde p\;,\quad \rho = \bar \rho + \tilde \rho\;.
\end{equation}
The fluid is assumed to be homogeneous and under the stationary flow described by $(\bar\wb,\bar p, \bar \rho)$ is considered as incompressible, so that a constant reference density $\bar \rho = \rho_0$ can be introduced. The acoustic perturbations are governed by the barotropic response denoting by $p_0$ and $\rho_0$ reference state variables, it holds that $p-p_0 = c^2(\rho - \rho_0)$, where $c^2$ is the squared acoustic velocity $c = \sqrt{k_f/\rho_0}$ with the bulk stiffness $k_f = 1/\gamma_f$, thus,  $\gamma_f$ is the fluid compressibility for the reference state. Upon substituting \eq{eq-NS2} in the Navier--Stokes equations, we obtain
\begin{equation}\label{eq-NS3}
  \rho_0  \left({\pdiff{}{t}}\tilde\ub + \bar\wb\cdot\nabla\tilde\ub + \tilde\ub\cdot\nabla\bar\wb\right) = - \nabla \tilde p + \nabla\cdot\tilde\taubf^\vis\;,\quad
   {\pdiff{}{t}}\tilde p + \wb\cdot\nabla \tilde p = - k_f \nabla\cdot\tilde\ub\;,
\end{equation}
while $(\bar\wb,\bar p)$ describing the stationary flow in the periodic scaffolds satisfy
\begin{equation}\label{eq-NS-flw}
  \rho_0\bar\wb\cdot\nabla\bar\wb + \nabla \bar p -\mu\nabla^2{\bar \wb} = \bar\fb\;,\quad
 \nabla\cdot \bar\wb = 0\;,
\end{equation}
where $\tilde\taubf^\vis$ is the viscous stress associated with the velocity fluctuations and $\bar\fb$ is the volume force.

In this paper, we consider further simplifying assumptions. Namely, the viscous effects are neglected, thereby  $\tilde\taubf^\vis \approx \zerobm$ in  \eq{eq-NS3}$_1$ and the steady flow \eq{eq-NS-flw} is treated as a potential flow, \ie $\bar\wb = -\nabla \Psi$   where the potential $\Psi$  satisfies $-\nabla^2 \Psi = 0$ in the fluid domain and $\nubf\cdot\nabla \Psi$ on the plate surface.
The issue of the flow problem homogenization is beyond the scope of the present paper, however, for the potential flow the homogenization procedure follows the one developed for the Helmholtz equation in the layer. In fact, the homogenized model can be derived relatively easily as  an exercise using the theoretical results reported in \cite{Rohan-Lukes-AMC2019}, or \cite{rohan-lukes-waves07}. 

Regarding the acoustic equations \eq{eq-NS3}, for stationary fluids, \ie when $\bar\wb \equiv 0$, these can be converted in the wave equation governing the pressure $\tilde p$.
When steady advection is respected, $\bar\wb \not \equiv \zerobm{}$, it is not straightforward  to eliminate the velocity in order to obtain a pressure formulation. Let us apply the divergence operator in \eq{eq-NS3}$_1$; this yields,
\begin{equation}\label{eq-NSadv2}
  \rho_0  \left(\nabla\cdot\dot{\tilde\ub} + 2 \pd_k \bar w_i \pd_i \tilde u_k +
  \bar\wb\cdot\nabla (\nabla\cdot\tilde\ub) + \tilde\ub\cdot\nabla(\nabla\cdot\bar\wb)
  \right) = -\nabla^2 \tilde p\;,
\end{equation}
where the second \lhs term can be approximated, as suggested in \cite{Rohan-Cimrman-AMC2021},  $\pd_i \bar w_k \pd_k \tilde u_i \approx \bar\wb\cdot\nabla(\nabla\cdot\tilde\ub)$. By the consequence, \eq{eq-NS3} governing the acoustic waves is approximated by the following equation,
\begin{equation}\label{eq-NSadv3}
  \left({\pdiff{}{t}} + \tau \bar\wb\cdot\nabla\right)\left(\pdiff{}{t}{\tilde p} + \bar\wb\cdot\nabla{\tilde p}\right) = c_f^2 \nabla^2 \tilde p\;,
\end{equation}
where $\tau = 3$ and $c_f = \sqrt{k_f/\rho_0}$. Further we define $\theta = (1+\tau)/2$,
thus, $\theta = 2$. We prefer to keep the abstract notation $\tau$ and $\theta$ in what follows. In the frequency domain, with $\FT{p}$ being considered as the Fourier image of $\tilde p$, \eq{eq-NSadv3} becomes an extension of the Helmholtz equation,
\begin{equation}\label{eq-Hlz}
  c_f^2 \nabla^2 \FT{p} + \om^2 \FT{p} - \imu\om(1+\tau) \pd_w \FT{p} - \tau \pd_w^2\FT{p} = 0\;,\quad \pd_w \FT{p} = \bar\wb\cdot\nabla\FT{p}\;,\quad \pd_w^2\FT{p} = \pd_w(\pd_w \FT{p})\;,
\end{equation}
where $\om$ is the frequency of incident waves, $\imu = \sqrt{-1}$.
The standard Helmholtz equation is obtained for $\bar\wb \equiv \zerobm$. Below, rather than $\FT{p}$ we simply use $p$ to refer to the acoustic fluctuations in the frequency domain. The notation $\pd_w$ and $\pd_w^2$ for the advective derivatives are employed throughout the next sections.

\subsection{Geometrical decomposition and the layer subproblem}

Although we are interested in modelling the acoustic field in a waveguide represented by a domain $\Om^G\subset \RR^3$ in which the perforated plate is embedded, we shall consider a  decomposition of  the ``global problem'' into two subproblems: the acoustic
interaction in the layer $\Om_\dlt$ and the outer acoustic problems in
$\Om^G\setminus \Om_\dlt$, where $\delta = \vkappa \veps$, with $\vkappa>0$ being fixed, is the layer thickness while $\veps$ characterizes the size of the plate perforations.
 The two subproblems are coupled by natural
 transmission conditions on the ``fictitious'' interfaces $\Gamma_\dlt^\pm$, see Fig.~\ref{fig:perflayer-scheme}.
The asymptotic analysis $\dlt\approx \veps\rightarrow 0$ is considered for the problem 
in $\Om_\dlt$ with the Neumann type boundary conditions on $\Gamma_\dlt^\pm$.  As the result of the dimensional reduction ``3D-to-2D'', the homogenized acoustic transmission layer is transformed in a problem defined on
$\Gamma_0$.
The  global problem for the acoustic
waves in the fluid interacting with the homogenized perforated plate
represented by $\Gamma_0$ is completed by a coupling condition which enables to introduce an implicit Dirichlet-to-Neumann operator.

\subsubsection{Decomposition of the global problem}\label{sec-decomp}

\begin{figure}[ht]
  \centering
  \includegraphics[width=0.98\linewidth]{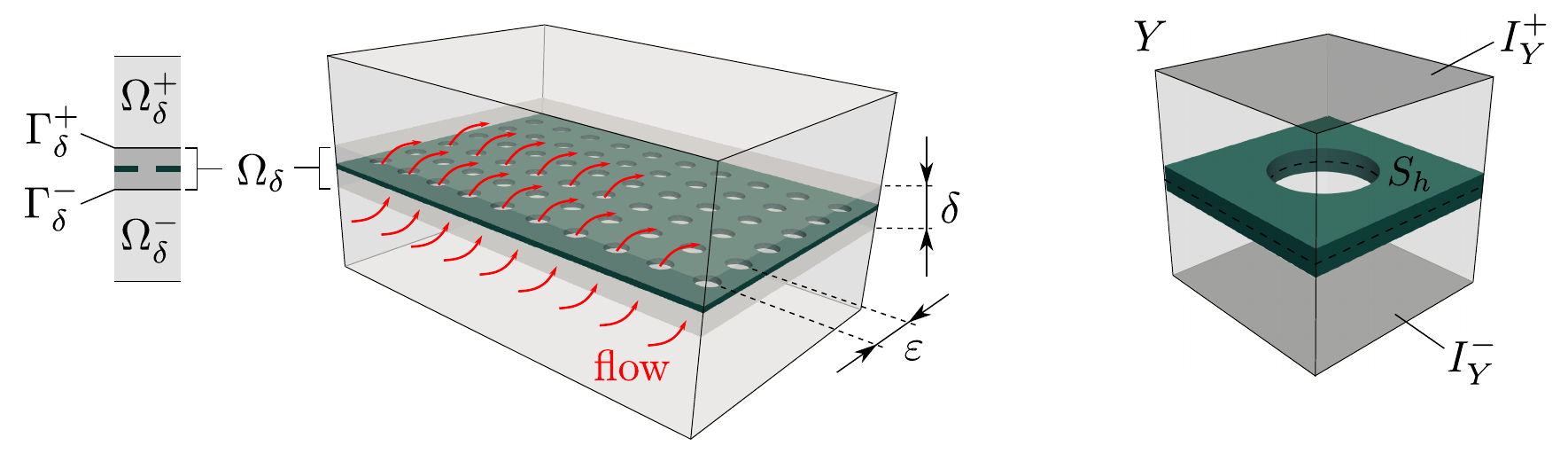}
  \caption{Left -- a waveguide and its domain decomposition in the layer $\Om_\dlt$ and
           the outer parts $\Om_\dlt^+$ and $\Om_\dlt^-$;
           right -- the representative cell.}  \label{fig:perflayer-scheme}
\end{figure}

Let $\Gamma_0 \subset \RR^2$ be an open bounded  domain 
spanned by coordinates $x_\alpha$, $\alpha = 1,2$, constituting a planar manifold in $\RR^3$. 
 Further let $\Gamma_\delta^{+}$ and
$\Gamma_\delta^{-}$ be equidistant to $\Gamma_0$ with the distance
$\delta/2 = {\rm dist}(\Gamma_0,\Gamma_\delta^+) = {\rm dist}(\Gamma_0,\Gamma_\delta^-)$.
We introduce $\Om_\delta = \Gamma_0 \times ]-\delta/2,\delta/2[ \subset \RR^3$, an open domain representing the transmission layer
bounded by $\pd \Om_\delta$ which is split as follows
\begin{equation}\label{eq-101}
\pd \Om_\delta = \Gamma_\delta^+ \cup  \Gamma_\delta^- \cup \pd_\ext \Om_\delta\;,
\quad \Gamma_\delta^\pm = \Gamma_0 \pm \frac{\delta}{2} \vec{e_3}\;,
\quad \pd_\ext \Om_\delta = \pd \Gamma_0 \times  ]-\delta/2,\delta/2[\;,
\end{equation}
where $\delta > 0$ is the layer thickness and $\vec{e_3} = (0,0,1)$,
see Fig.~\ref{fig:perflayer-scheme}.

In the waveguide $\Om^G$, the fluid occupies domain  $\Om_\dlt^{G,\veps} = \Om_\dlt^+ \cup \Om_\dlt^-\cup \Om^{*\veps}$, where $\Om^{*\veps}\subset\Om_\dlt$ is the fluid-saturated  part of the  the ``transmission layer'' $\Om_\dlt$  containing the fluid and the rigid perforated plate $\Sigma^\veps \subset \Om_\dlt$.
According to the decomposition, the acoustic fields 
represented by $P^\dlt$ in $\Om_\dlt^\pm$ and $p^\veps$ in $\Om^{*\veps}$
both satisfy the extended Helmholtz equation \eq{eq-Hlz}, being coupled on $\Gamma_\dlt^\pm$ by
\begin{equation}\label{eq-G2-01}
    \pdnw{P^\dlt} = \imu\om g^{\veps\pm} \quad \mbox{ on } \Gamma_\dlt^\pm\;, \qquad
P^\dlt = p^\veps \quad \mbox{ on } \Gamma_\dlt^\pm\;,\\
\end{equation}
where $\pdnw{}$ is the advection normal derivative depending on $\om$ and $\bar \wb$,
\begin{equation}\label{eq-pdnw}
    \pdnw{q} = \pd_n q - \bar w_n (\imu\om \theta q + \tau \pd_w q) \;,
\end{equation}
involving projections $\pd_n = \nb\cdot\nabla$, $w_n = \nb\cdot\wb$, 
and $g^{\veps\pm}$ represents the acoustic momentum flux involved in the definition of the layer subproblem \eq{eq-G2-02}.
The definition of $\pdnw{}$ is obtained when deriving the weak formulation of the acoustic problem featured by the fluid advection.
Without loss of generality, when dealing with the unit normal vector $\nb$ on $\Gamma_\dlt^\pm$, we consider its orientation outward to layer $\Om_\dlt$.



The acoustic potential $p^\veps$ in the layer satisfies 
\begin{equation}\label{eq-G2-02}
  \begin{split}
    c_f^2 \nabla^2 p^\veps + \om^2 p^\veps - \imu\om(1+\tau) \pd_w p^\veps - \tau \pd_w^2 p^\veps  & = 0 \quad \mbox{ in } \Om_\dlt^\veps \;, \\
\mbox{  \textbf{interface conditions}  } \quad \pdnw{p^\veps} & = - \imu \om g^{\veps\pm} \quad \mbox{ on } \Gamma_\dlt^{\pm}\;,\\
p^\veps & = P^\dlt  \quad \mbox{ on } \Gamma_\dlt^\pm\;,\\
\mbox{ zero velocity of solid structure,  } \quad \pdiff{p^\veps}{n} & = 0  \quad \mbox{ on } \pd \Sigma^\veps\;, \\
\pdiff{p^\veps}{n} & = 0  \quad \mbox{ on } \pd_\ext\Om_{\delta}\;,
\end{split}
\end{equation}
where $c$ is the speed of sound propagation. Note that $\pdiff{}{n} = \nb\cdot\nabla = \pdnw{}$ on the solid walls, since there $\bar w_n = 0$. 

Although we do not specify the outer problem in $\Om_\dlt^+ \cup \Om_\dlt^-$, the advection normal derivative must be taken into account when dealing  with the open waveguide boundary conditions (incident waves, or non-reflective boundaries).

\subsubsection{Periodic microstructure in the layer}\label{sec-geometry}

The microstructure of the perforation is periodic, being generated by
the representative  periodic cell  $Y_\vkappa = \Xi \times \vkappa {]-1/2,+1/2]}$, where
$\Xi = {]0,b_1[}\times{]0,b_2[}$. We can assume $|\Xi| = b_1b_2 \approx 1$, hence $|Y_\vkappa| \approx \vkappa$. The surface of the hexahedron $Y_\vkappa$ is decomposed into
the upper and lower sides, $I_Y^+$ and $I_Y^-$, and the ``periodicity'' sides $\pd_\# Y_\vkappa$, thus $\pd Y_\vkappa =I_Y^+ \cup  I_Y^- \cup \pd_\# Y_\vkappa$.

The solid (plate) is introduced by a representative solid obstacle $S_h \subset  \Xi \times h{]-1/2,+1/2]}$, where $h << \vkappa$ is the plate thickness, \ie the maximum thickness of the obstacle (measured transversally to the layer thickness). 

The fluid domain $\Om^{*\veps}$ is generated  using $Y_\vkappa^* = Y_\vkappa \setminus \ol{S}_h$ as a $\Xi$-periodic lattice,
\begin{equation}\label{eq-def-layer}
\Om^{*\veps} = \bigcup_{\zeta' \in \ZZ^2}\veps(\zeta + Y_\vkappa^*)\cap \Om_\delta\;.\quad \zeta = (\zeta',0)\;.
\end{equation}

In what follows, we refer to $\Xi$-periodic functions: any such function, say $f(x,y)$ where $x\in\Om$, $y\in Y$, satisfies $f(\cdot,y) = f(\cdot,y + w)$ for $w = k_1 b_1 \vec e_1 + k_2 b_2 \vec e_2$ with $k_1,k_2 \in \ZZ$. By $H_\plper^1(Y) \subset H^1(Y)$ we denote the subspace of $\Xi$-periodic functions.
We shall also simplify the notation $Y\equiv Y_\vkappa$, since $\vkappa$ is a fixed constant.

\section{The layer homogenization}\label{sec-homog}
In this section, we introduce the limit two-scale equations of the acoustic problem imposed in the
transmission layer. We employed the asymptotic analysis based on the unfolding method \cite{Cioranescu2008-Neumann-sieve} which uses the standard convergence results. For this, a~priori estimations on the acoustic pressure are needed; theses can be obtained in analogy with the treatment of the vibroacoustic problem \cite{Rohan-Lukes-AMC2019}. 
In our setting, the unfolding operator $\Tuftxt{}:  L^2(\Om_\delta;\RR) \rightarrow L^2(\Gamma_0 \times Y;\RR)$ transforms a function $f(x)$ defined in $\Om_\delta$ into a function of two variables, $x' \in \Gamma_0$ and $y \in Y$. 
For any $f \in L^1(Y)$, the cell average involved in all unfolding integration formulae will be abbreviated by
\begin{equation}\label{eq-uf4}
\frac{1}{|\Xi|}\int_\Xi f= \intY_\Xi f\;,\quad \frac{1}{|\Xi|}\int_{D} f =:  \intY_D f\;,
\end{equation}
whatever the domain $D \subset \ol{Y}$ of the the integral is (\ie volume, or surface). In what follows, all variables depending on the scale $\veps$ are labelled by $^\veps$.

In contrast with previous studies of the acoustic transmission, where a static fluid was considered, the flow phenomenon not only modifies the acoustic wave model, but also requires to account for a locally periodic structures.
Independently of the flow model used to analyze and compute the advection velocity filed $\wb^\veps$, we assume that such a model provides a bounded unfolded velocity field, such that for a real $\beta \geq 0$
\begin{equation}\label{eq-ufw}
  \veps^\beta\Tuf{\wb^\veps} \cwto \wb(x',y) \quad \mbox{ weakly in } L^2(\Gamma_0 \times Y^*)\;, \quad 
  \wb \mbox{ is bounded in  } L^\infty(\Gamma_0 \times Y^*)\;.
\end{equation}
This assumption on the advection velocity ensures that convergence results reported below can be obtained in much similar way as in the cases of a static fluid when the standard Helmholtz equation is employed instead of \eq{eq-Hlz}.
However, by the consequence, the ``microconfigurations'' are only locally periodic being featured by $\wb(x',\cdot)$ depending on $x'\in \Gamma_0$ and constituting the differential operator.

\subsection{Variational formulation and the convergence results} 
When deriving the weak formulation  of the layer subproblem \eq{eq-G2-02}, the advection derivative \eq{eq-pdnw} is obtained, being substituted by the interface coupling condition  involving $g^{\veps\pm}$. The acoustic pressure $p^\veps \in  H^1(\Om^{*\veps})$ is a weak solution of \eq{eq-G2-02}, iff
\begin{equation}\label{eq-va6F}
\begin{split}
& c^2 \int_{\Om^{*\veps}} \nabla p^\veps \cdot \nabla q^\veps
  - \om^2 \int_{\Om^{*\veps}} p^\veps q^\veps
  \blue{+ \imu\om\theta\int_{\Om^{*\veps}}(q^\veps \pd_w p^\veps - \pd_w q^\veps p^\veps) - \tau\int_{\Om^{*\veps}}\pd_w q^\veps  \pd_w p^\veps}\\
  & =
-\imu \om c^2\left ( \int_{\Gamma^{\pm\veps}} \magenta{g^{\veps\pm}} q^\veps\,d\Gamma
 \right ),\quad \mbox{ for all } q^\veps \in H^1(\Om^{*\veps})\;.
\end{split}
\end{equation}
The \rhs terms depending on the acoustic momentum fluxes must be specified. Let $g^{0} \in L^2(\Gamma_0)$ and
$g^{1\pm}(x',y') \in L^2(\Gamma_0 \times \RR^2)$, whereby $g^{1\pm}(x',\cdot)$ being $\Xi$-periodic in the second
variable. For 
\begin{equation}\label{eq-aes2}
\hat g^{\veps+}(x') = g^{0}(x')  + \veps g^{1+}(x',\frac{x'}{\veps}) \;, \quad
\hat g^{\veps-}(x') = -g^{0}(x')  - \veps g^{1-}(x',\frac{x'}{\veps})\;,
\end{equation}
the following convergence results hold
\begin{equation}\label{eq-aes5a}
\begin{split}
p^\veps & \cwto p^0 \quad \mbox{ weakly in } L^2(\hat\Om)\;,\quad
\pd_z p^\veps \cwto 0\quad \mbox{ weakly  in } L^2(\hat\Om)\;,
\end{split}
\end{equation}
thus, $\pd_z p^0 = 0$. Moreover, the classical results of the homogenization yield
\begin{equation}\label{eq-aes5}
\begin{split}
\Tuf{p^\veps} & \cwto p^0 \quad \mbox{ weakly in } L^2(\Gamma_0 \times Y^*)\;,\\
\Tuf{\gradpl p^\veps} & \cwto \gradpl_{x'} p^0 + \gradpl_{y'}p^1 \quad \mbox{ weakly in } L^2(\Gamma_0 \times Y^*)\;,\\
\frac{1}{\veps}\Tuf{\pd_z p^\veps} & \cwto  \pd_z p^1 \quad \mbox{ weakly in } L^2(\Gamma_0 \times Y^*)\;.
\end{split}
\end{equation}

\subsection{Limit two-scale equations}\label{sec:two-scale-equations}
These can be obtained rigorously using the convergence results \eq{eq-aes5a}-\eq{eq-aes5} enabling the asymptotic analysis to be applied directly in \eq{eq-va6F}. Formal procedure simplifying the
derivation of limit equations is based on the matched asymptotic expansions which require the use of the recovery sequences constructed in accordance
with \eq{eq-aes5a}-\eq{eq-aes5}. Neglecting the higher order terms in $\veps$, we consider  the
following approximate expansions for unfolded unknown $p^\veps$ and the test function $q^\veps$,
\begin{equation}\label{eq-va13a}
\Tuf{p^\veps} = p^0(x') + \veps p^1(x',y)\;, \quad \Tuf{q^\veps}  = q^0(x') + \veps q^1(x',y)\;, \\
\end{equation}
where $x' \in \Gamma_0$, $ y' \in \Xi$ and $y = (y',z) \in Y$; in~\eq{eq-va13a}, recalling that 
all the two-scale functions are $\Xi$-periodic in the second variable.
The detailed derivation of the limit problem equations is skipped in this paper, since the resulting equations can be obtained in analogy with the treatment reported in \cite{Rohan-Lukes-AMC2019}, or \cite{Rohan-Lukes-AMM2022}.

For this technical procedure and the result presentation the following notation is needed,
\begin{equation}\label{eq-w-grad}
\begin{split}
  \gradplx p = (\pd_\alpha^x p) = (\pd_1^x,\pd_2^x) p\\
 \bar\pd_w^x p = \wb\cdot\gradplx p = \pl{\wb}\cdot\gradplx p = (w_1\pd_1^x,w_2\pd_2^x) p\\
 \bar\pd_w^y p = \wb\cdot\gradply p = \pl{\wb}\cdot\gradply p = (w_1\pd_1^y,w_2\pd_2^y) p\\
  \pd_w^y p = \wb\cdot\nabla_y p = (w_1\pd_1^y,w_2\pd_2^y,w_3\pd_3^y) p = (w_1\pd_1^y,w_2\pd_2^y,w_3\pd_z) p \\
  \pd_w^y y_\alpha = (\dlt_{1\alpha},\dlt_{2\alpha},0)\\
  \pd_w^y \pi^\alpha = \wb\cdot\nabla_y \pi^\alpha\;,
\end{split}
\end{equation}
where $\wb = (\wb',w_3) = (w_1,w_2,w_3)$. In the expressions presented below, a product of $\ab\in\RR^3$ and $\bb\in\RR^2$ is correct $\ab\cdot\bb = \sum_\alpha a_\alpha b_\alpha$ when $a_3 = 0$, which justifies the writing in \eq{eq-w-grad}$_{2,3}$.

The two-scale asymptotic homogenization method applied to \eq{eq-va6F} yields the limit two-scale variational equation 
\begin{equation}\label{eq-f6}
\begin{split}
& c^2 \int_{\Gamma_0} \intY_{Y^*} (\gradplx p^0 + \gradply p^1)\cdot
  (\gradplx q^0 + \gradply q^1) + c^2 \int_{\Gamma_0} \intY_{Y^*}\pd_z p^1 \pd_z q^1 - \om^2 \int_{\Gamma_0} \intY_{Y^*} p^0 q^0\\
 & \dgreen{-\tau \int_{\Gamma_0} \intY_{Y^*}(\bar\pd_w^x p^0 + \pd_w^y p^1)(\bar\pd_w^x q^0 + \pd_w^y q^1)}\\
  & \dgreen{ + \imu\om\theta  \int_{\Gamma_0} \intY_{Y^*} \left(
  q^0(\bar\pd_w^x p^0 + \pd_w^y p^1) - p^0(\bar\pd_w^x q^0 + \pd_w^y q^1) \right) }\\
&  = -\imu\om c^2 \int_{\Gamma_0} \left[
q^0 \intY_\Xi \magenta{\Delta g^{1}} + \red{g^{0}}\left(\intY_{I_y^+}q^1 -  \intY_{I_y^-}q^1\right)\right]\\
\end{split}
\end{equation}
for all $q^0 \in V_0(\Gamma_0),q^1\in H_\plper^1(Y^*)$. Thus, the fluid flow modifies the limit problem by the two integrals involving the advection derivatives.

The limit model of the layer is coupled with the outer acoustic fields $P^\dlt$ due to conditions \eq{eq-G2-01}; to respect  \eq{eq-G2-01}$_2$, we consider its weak form, see the identity (61) in \cite{Rohan-Lukes-AMC2019}. Its approximation for a given finite heterogeneity scale $\veps_0$, which is related to a finite layer thickness
$\delta_0 = \vkappa\veps_0 >0$, yields the following limit condition
\begin{equation}\label{eq-cc3}
\frac{1}{\veps_0} \int_{\Gamma_0} \psi (\hat P^{+} - \hat P^{-} ) = \int_{\Gamma_0} \psi \left(\intY_{I_y^+} p^1 - \intY_{I_y^-} p^1\right)
 \quad \forall \psi \in L^2(\Gamma_0)\;,
\end{equation}
where $\hat P^{+/-}$ are the limit traces of $P^\dlt$ on $\Gamma_\dlt^{+/-}$ for $\dlt \rightarrow 0$.




\subsubsection{Autonomous problems} 

When testing the limit equation \eq{eq-f6} with $q^1\not=0$ while $q^0 = 0$,
the local problem in the fluid part is obtained. It provides characteristic responses and thereby also the homogenized coefficients constituting the macroscopic model of the homogenized layer.
$p^1$ depends on the ``macroscopic'' functions $\pd_\alpha^x p^0,p^0$ and $g^0$. Therefore,  due to the linearity, we
can introduce the following split
\begin{equation}\label{eq-f7}
p^1(x',y) =  \blue{\pi^\beta(y)}\pd_\beta^x p^0(x') + \imu\om  \left(\blue{\xi(y)} g^0(x') + \dgreen{\pi^P(y)} p^0\right)\;.\\
\end{equation}

Define the {operator $\AA$} and the inner product for any functions \blue{$p,q \in H_\plper^1(Y^*)$},
\begin{equation}\label{eq-lp2}
  \ipYf{\AA p}{q} := \blue{\ipYs{\nabla_y p}{\nabla_y q}}
-\frac{\tau}{c^2} \dgreen{\ipYs{\pd_w^y  p}{ \pd_w^y q}}
  \;.
\end{equation}
Due to the involvement of $\wb(x',\cdot)$ in \eq{eq-lp2},  $\AA$ is parameterized by $x' \in \Gamma_0$.
Using the split form substituted in \eq{eq-f6}, 
three local problems are distinguished. 

Their solutions provide local characteristic responses of the layer microstructure \wrt macroscopic quantities $\gradplx p^0, g^0$ and $p^0$
\begin{itemize}
\item Find $\pi^\beta \in H_\plper^1(Y^*)$, such that \\ 
\begin{equation}\label{eq-lp3}
     \ipYf{\AA\pi^\beta}{q} = - \ipYs{\gradply y_\beta}{\gradply q}
    +\dgreen{\frac{\tau}{c^2}\ipYs{\pd_w^y y_\beta}{\pd_w^y q}}\;, \forall q \in  H_\plper^1(Y^*)\;.
\end{equation}
\item Find $\xi \in H_\plper^1(Y^*)$, such that \\ 
\begin{equation}\label{eq-lp4}
    \ipYf{\AA\xi}{q} = -\left(\intY_{I_y^+} q
- \intY_{I_y^-} q\right) \;,\quad \forall q \in H_\plper^1(Y^*)\;.
\end{equation}

\item Find $\pi^P \in H_\plper^1(Y^*)$, such that  \\ 
\begin{equation}\label{eq-lp5}
  \ipYf{\AA\pi^P}{q} = \frac{\theta}{c_f^2}\int_{Y^*} \wb\cdot \nabla_y q \;,\quad \forall q \in H_\plper^1(Y^*) \;.    
\end{equation}
\end{itemize}

\begin{myremark}{rem1}
  Note that problem \eq{eq-lp5} appears only for nonvanishing flow, \ie if $\wb \not \equiv \zerobm$, while the other two problems are relevant also in the static fluid acoustics.
  Since, in general, $\wb(x',y)$ is merely locally periodic, all the three local problems are specific to a vicinity of the macroscopic position $x' \in \Gamma_0$.
  \end{myremark}

\subsection{Macroscopic transmission layer model}\label{sec:macro-layer-model}
Upon substituting the split forms \eq{eq-f7} in the limit equation \eq{eq-f6}, the macroscopic equation expressed in terms of the homogenized coefficients $\Ab,\Bb,M_w,T_w, \bar\Wb'$, and $\bar\Wb$ can be obtained:
\begin{equation} \label{eq-M1}
\begin{split}
& c^2\int_{\Gamma_0} ( \blue{\Ab} \gradpl_x p^0 )\cdot \gradpl_x q^0
  - \om^2(\zeta^* + \magenta{M_w})  \int_{\Gamma_0} p^0 q^0 + \imu \om c^2\int_{\Gamma_0} g^0 \blue{\Bb}\cdot\gradpl_x q^0 \\
& - \om^2\theta \int_{\Gamma_0} q^0 \magenta{T_w} g^0 
  + \imu \om \theta\int_{\Gamma_0}\left(q^0 \magenta{\bar\Wb }\cdot\gradpl_x p^0 + \gradplx q^0 \blue{\bar\Wb'} )p^0\right)
+ \imu \om c^2 \int_{\Gamma_0} q^0 \Dlt G^1  = 0\;,
\end{split}
\end{equation}
for all $q^0 \in H^1(\Gamma_0)$, where $\zeta^* = |Y^*|/|Y|$, and $\Delta G^{1}:=\intY_{\Xi}(g^{1+}-g^{1-})$.
In addition, the coupling equation involving further coefficients $\Bb',T_w'$ and $F$ is derived in the similar way from \eq{eq-cc3},
\begin{equation}\label{eq-M2}
\int_{\Gamma_0} \psi \left(
\Bb'\cdot\gradpl_x p^0 - \imu \om F g^0 + \imu \om T_w' p^0
\right) = \frac{1}{\veps_0} \int_{\Gamma_0} \psi (\Delta{ P^{\dlt_0}})\quad \forall \psi \in L^2(\Gamma_0)\;.
\end{equation}
where $\Delta{P^{\dlt_0}} = P^{\dlt_0+} -  P^{\dlt_0-}$ is the jump of the solution in the outer part $\Om_{\dlt_0}$. In \Appx{sec-appx1}, we give expressions for all the homogenized coefficients satisfying the following symmetry relationships,
\begin{equation}\label{eq-M3}
A_{\alpha\beta} = A_{\beta\alpha}\;,\quad B_\alpha = B_\alpha'\;, \quad T_w'=-\frac{\theta}{c^2}T_w\;, \quad \bar\Wb' = -\bar\Wb\;.
\end{equation}
This property leads to the Hermitean symmetry of the discretized system  \eq{eq-M1}-\eq{eq-M2} presented in the matrix form.
Besides the effective advection flow velocity $\bar\Wb$, also other coefficients labelled by subscript $w$ vanish for a static fluid.

\section{Numerical simulations}\label{sec:num-simulations}


In this section we illustrate the influence of the flowing fluid on acoustic waves propagating in a waveguide fitted with the perforated rigid plate. Recall that the acoustic problem \eq{eq-M1}-\eq{eq-M2} imposed in the layer involves variables $p^0,g^0$,
$\Delta G^1$ and $\Delta P^{\dlt_0}$. However, as explained in \cite{Rohan-Lukes-AMC2019},
$p^0$ and $g^0$ can be replaced by $G^+,G^-$ and $P^{\dlt_0+},P^{\dlt_0-}$, such that $p^0 = (P^{\dlt_0+}+P^{\dlt_0-})/2$, and $g^0 = (G^+ + G^-)/2$, whereas $\Dlt \hat P = P^{\dlt_0+}-P^{\dlt_0-}$ and $\Dlt G^1 = (G^+ - G^-)/\veps_0$ depending on scale $\veps_0 >0$ determining the layer thickness $\dlt_0 = \vkappa\veps_0$. The global model of the waveguide  consists of the following three parts:
\begin{list}{}{}
\item {\it (i)} The outer part model, eq. \eq{eq-Hlz} governing $P^{\dlt_0}$ in $\Om_{\dlt_0}^\pm := \Om_{\dlt_0}^+\cup \Om_{\dlt_0}^-$ complemented by boundary conditions on $\pd_\ext\Om_{\dlt_0}^\pm:=\pd\Om_{\dlt_0}^\pm\setminus\Gamma_{\dlt_0}^\pm$.
\item {\it (ii)} The homogenized layer model, \eq{eq-M1}-\eq{eq-M2}, where the above substitutions of  $p^0$ and $g^0$ are used.
 \item {\it (iii)} The interface conditions \eq{eq-G2-01}$_1$  applied at $\Gamma_{\dlt_0}^\pm$, where $g^{\veps\pm}$ is replaced by $g^0$.
\end{list}
These three coupled parts constitute the problem for $P^{\dlt_0}$ and $G^+,G^-$. Its
numerical solutions are computed by a monolithic approach using the finite element (FE) method implemented in the Python
based package {\it SfePy}: Simple Finite Elements in Python.
Note that the homogenized model \eq{eq-M1}-\eq{eq-M2} depend on the microscopic autonomous problems which are well decoupled from the global problem solution, as described in Sections~\ref{sec:two-scale-equations} and \ref{sec:macro-layer-model}. This is the main advantage of the homogenized interface model, since the cell problems are cheap to solve independently of the global solutions. Discretization of the finite scale layer involving the perforated plate is avoided. All the unknown fields defined in $\Om$, on $\Gamma_0$ and in $Y^*$ are approximated by linear Lagrangian finite elements.

\begin{figure}[ht]
    \centering
    \includegraphics[width=0.98\linewidth]{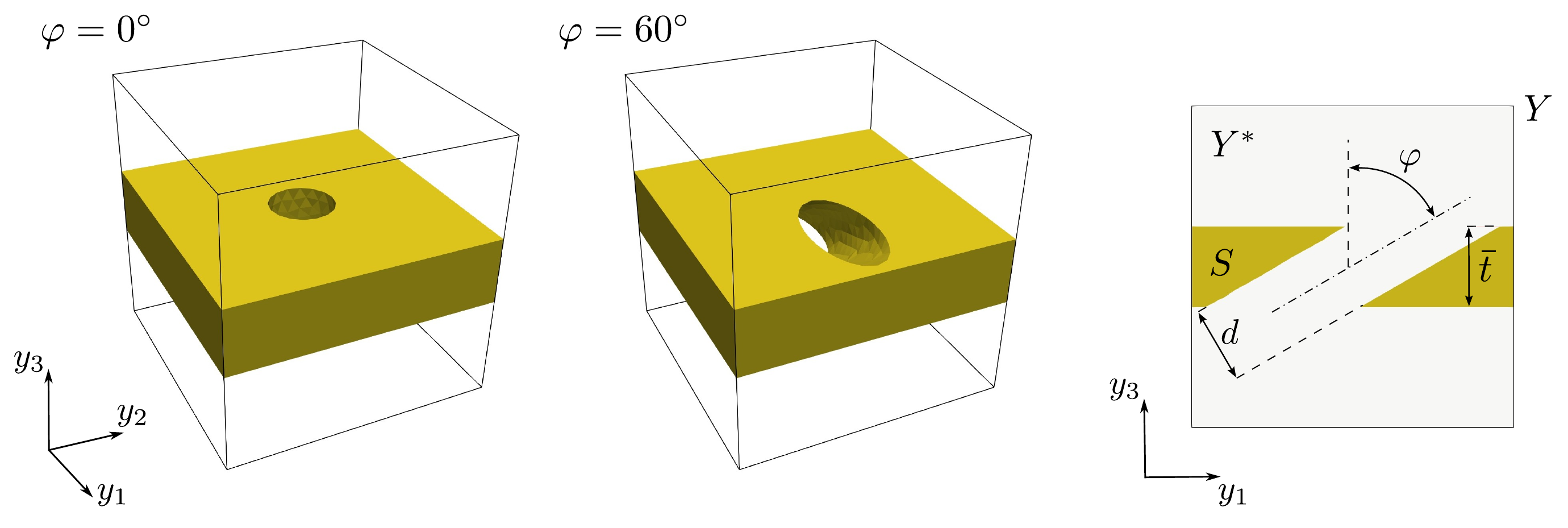}
    \caption{Periodic unit cell $Y$ -- a rigid plate perforated by a cylindrical hole declined by angle $\vphi$.}
    \label{fig:num-micro-geom}
\end{figure}

In the examples reported below, we consider plates perforated by cylindrical holes characterized by the diameter $d = 0.24\,\veps_0$\,m and the axis
slope $\varphi$, see Fig.~\ref{fig:num-micro-geom}. The acoustic fluid
occupying domain $Y^\ast$ is characterized by its density $\varrho_0 =
1.55\,\mbox{kg/m}^{3}$ and by the sound speed $c = 343\,\mbox{m/s}$.



\subsection{Homogenized acoustic coefficients -- dependence on the flow velocity}
\begin{figure}[ht]
    \centering
    \includegraphics[width=0.99\linewidth]{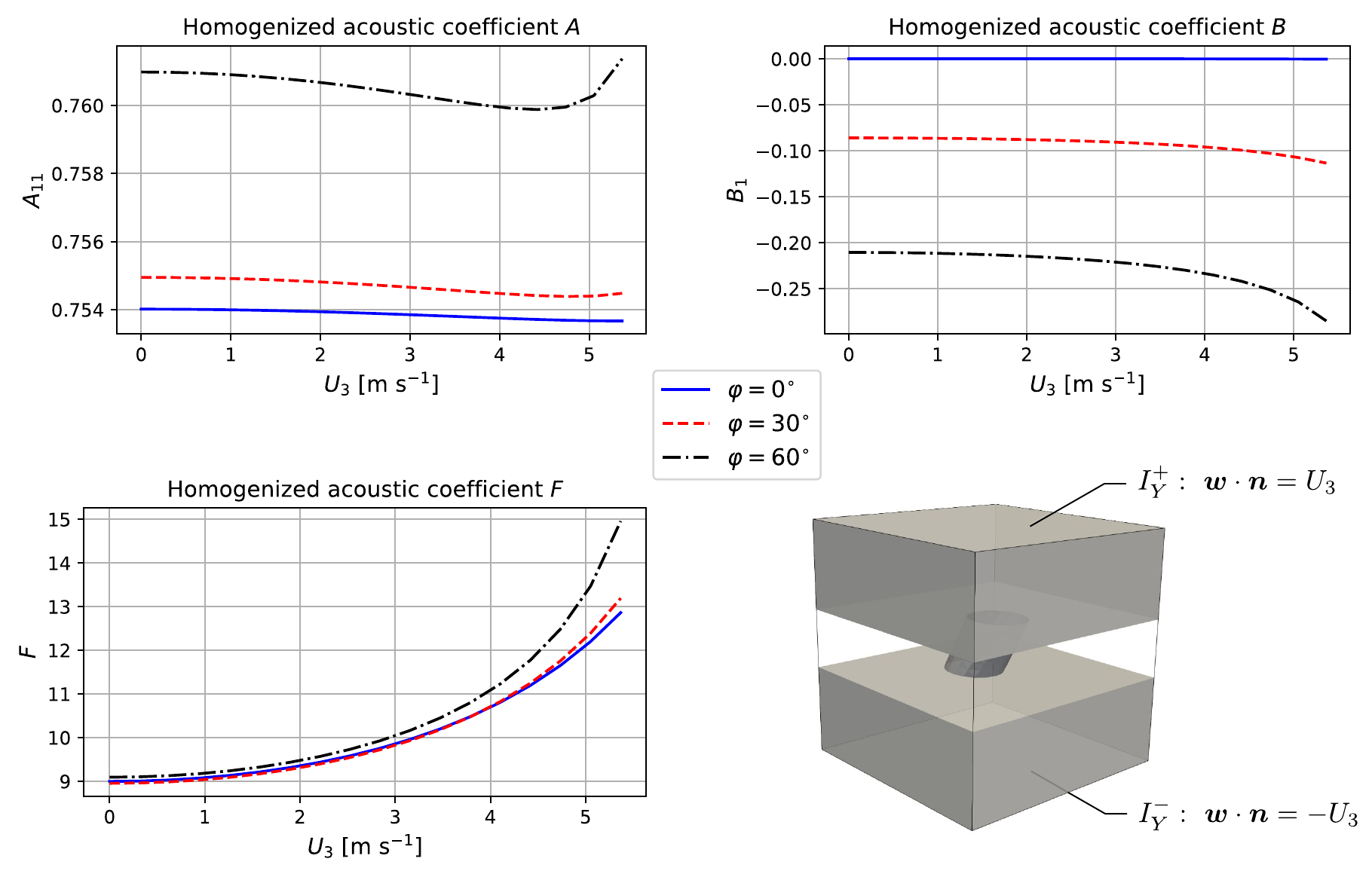}
    \caption{Dependence of the homogenized coefficients $A_{11}$, $B_1$, and $F$
        on the prescribed boundary velocity $U_3$ and
        on the geometrical parameter $\varphi$.}
    \label{fig:num-coefs-cv-phi}
\end{figure}

In addition to the geometry of the perforations, as represented by domain
$Y^\ast$, the homogenized acoustic coefficients  \eq{eq-F5}--\eq{eq-F7b} of the macroscopic equations \eq{eq-M1}-\eq{eq-M2} also depend on the flow
velocity which is involved in the local problems \eq{eq-lp3}--\eq{eq-lp5}. Fig.~\ref{fig:num-coefs-cv-phi}
shows the dependence of the coefficients on the
velocity $\wb$ and on the angle parametrizing the hole slope,
$\vphi = 0^{\circ}, 30^{\circ}, 60^{\circ}$, see Fig.~\ref{fig:num-micro-geom}. At the microlevel, velocity field $\wb$ is obtained due to the reconstruction of macroscopic flow velocity, in general; this issue is beyond the scope of the present paper, although, in the example of the global acoustic filed in a waveguide, the homogenized model of the potential flow is employed. However, in this study, it is sufficient to solve an auxiliary potential flow problem for a $\Xi$-periodic velocity $\wb \in H_\plper^1(Y^*)$ which is a weak solution of
$$
\wb = -\nabla\Phi, \qquad \nabla \cdot \wb = 0\quad \mbox{ in }Y^*\;, \quad \wb \cdot \nb = \pm  U_3\quad \mbox{ on }I_Y^\pm
$$
%
with a prescribed velocity $U_3$ varying in range from $0\,\mbox{m/s}$ to $5.5\,\mbox{m/s}$.

\begin{figure}[ht]
    \centering
    \includegraphics[width=0.8\linewidth]{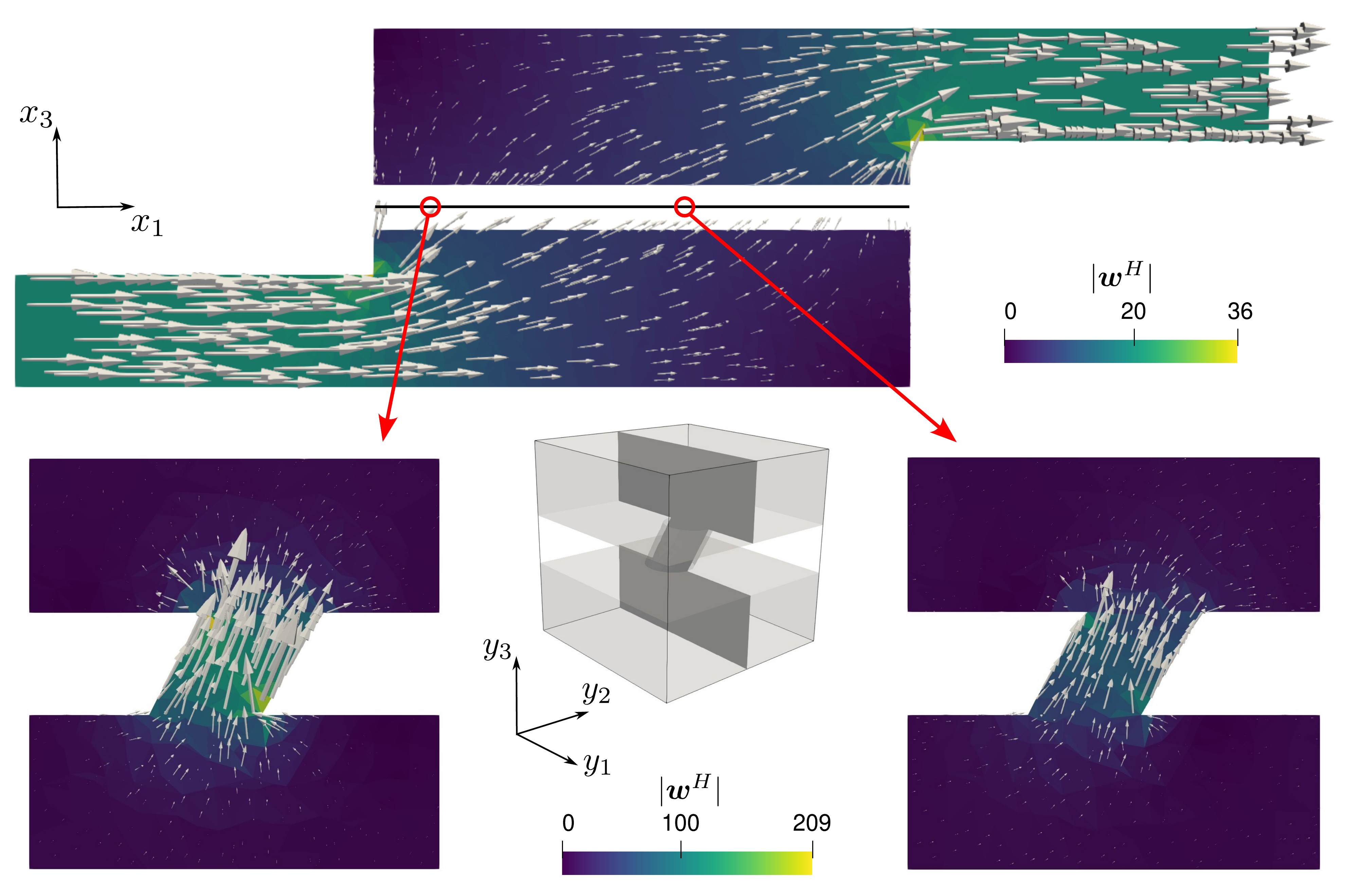}
    \caption{Potential flow: top -- the macroscopic flow velocity;
                             bottom -- the reconstructed velocity field at the microscopic level
                                       in two distinct macroscopic points.}
    \label{fig:num-flow-velocity}
\end{figure}

\subsection{Wave propagation in a waveguide with the homogenized interface layer}\label{sec:two-scale-problem}

\begin{figure}[ht]
    \centering
    \includegraphics[width=0.98\linewidth]{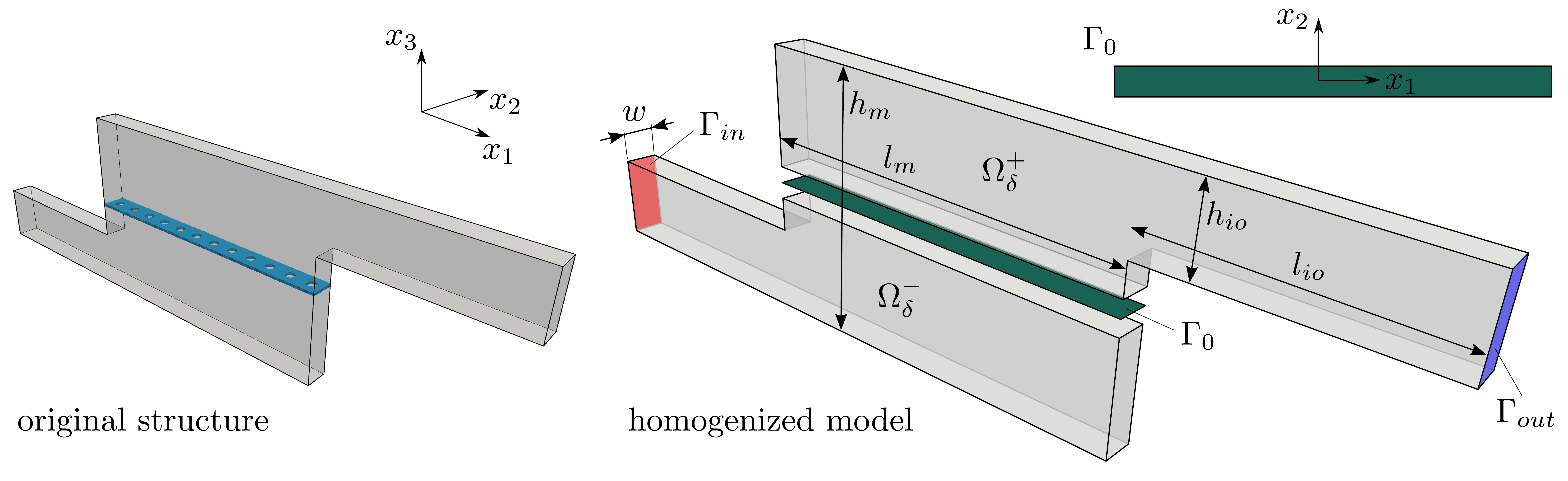}
    \caption{Waveguide -- decomposition of the macroscopic domain.}
    \label{fig:num-macro-geom}
\end{figure}

We consider an acoustic waveguide consisting of two equally shaped parts,
separated by the perforated plate. According to the global problem
decomposition, as announced in Section~\ref{sec-decomp}, the acoustic field
$P^{\dlt_0}$ is defined in domains $\Om_{\dlt_0}^+$ and $\Om_{\dlt_0}^-$,
whereas the transmission layer $\Om_{\dlt_0}$ of the thickness $\dlt_0 =
\veps_0 = 0.025$ separating the two domains is represented by the homogenized
interface $\Gamma_0$. see Fig.~\ref{fig:num-macro-geom}. The waveguide is
characterized by dimensions $l_m = 0.3$\,m, $h_m = l_{io} = 0.2$\,m, $h_{io} =
0.0625$\,m, and $w = 0.01$\,m. For simplicity, in the $x_2$-direction, we apply
the periodic boundary conditions on the faces orthogonal to the $x_2$-axis,
which yields a homogeneous distribution of the macroscopic fields \wrt $x_2$.
Hence the macroscopic problem is quasi 2D, which also facilitate the
visualization and interpretation of obtained results.

\begin{figure}[ht]
    \centering
    \includegraphics[width=0.8\linewidth]{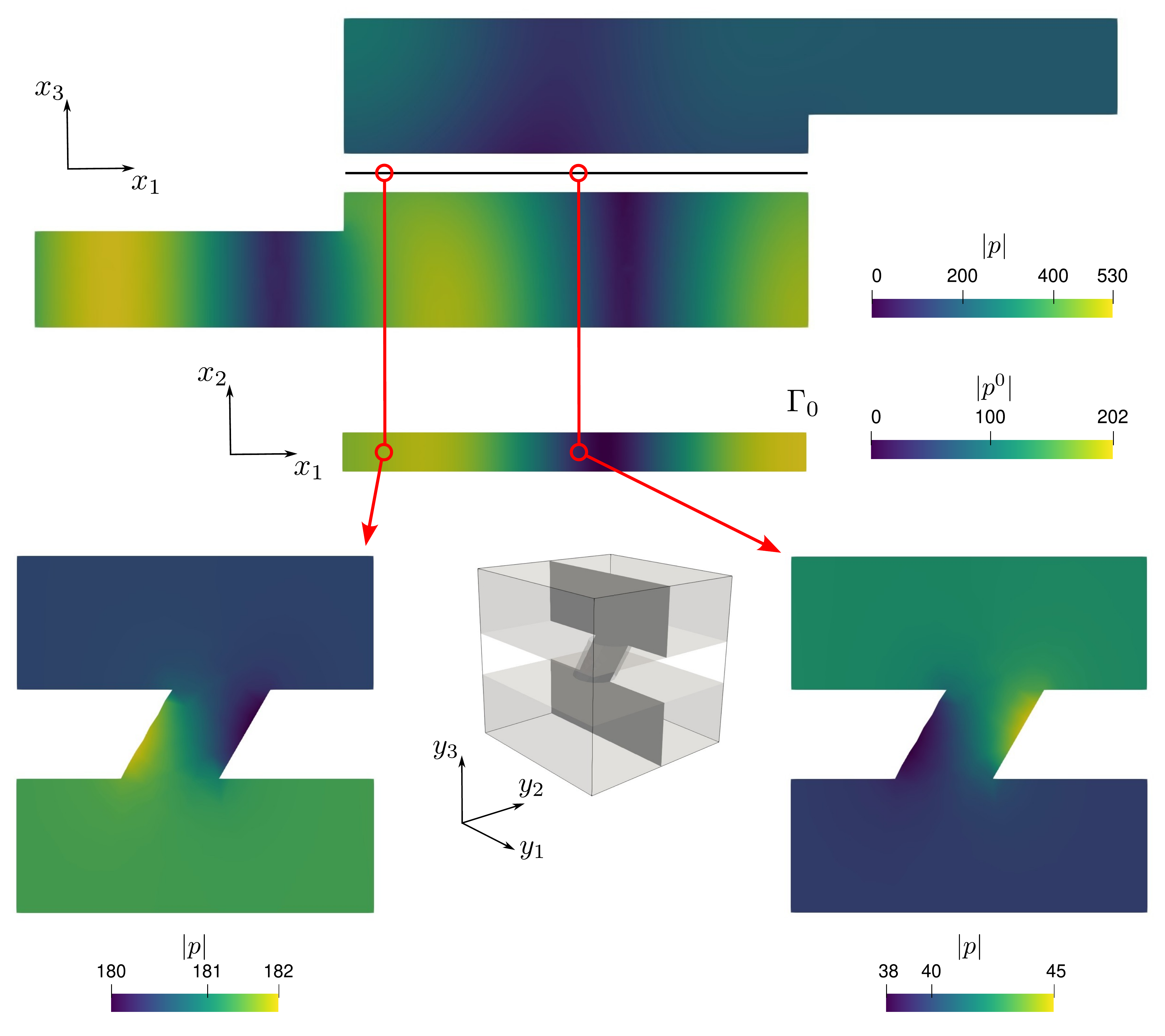}
    \caption{Magnitude of the acoustic pressure:
        top -- pressure in section of $\Omega$;
        middle -- acoustic pressure $p^0$ in $\Gamma_0$;
        bottom -- the pressure fields reconstructed at the microscopic level
                                       in two distinct macroscopic points.}
    \label{fig:num-pressure}
\end{figure}

The periodic interface is represented by a periodic unit cell in which the
rigid plate (domain $S$) of thickness $\bar{t} = 0.25\,\veps_0$\,m serves an
obstacle for the flow and the acoustic waves. The scaling parameter is chosen
as $\varepsilon = 0.3 / 12$ which corresponds to a perforated interface formed
of 12 periodic cells in the $x_1$-direction. In
Fig.~\ref{fig:num-flow-velocity}, we display a solution of the potential flow
problem imposed in the waveguide for the uniform inlet velocity $U_{in} =
-U_{out} = 20\,\mbox{m/s}$ on $\Gamma_{in}$. This was computed using a
two-scale model derived by the homogenization of the transmission layer
$\Om_\dlt$ using an analogical procedure to the one employed in the acoustic
problem homogenization. By virtue of the unfolding procedure, the local flow
field at the microlevel can be reconstructed. Using the flow velocities so
obtained we solve the microscopic subproblems \eq{eq-lp3}--\eq{eq-lp4} and
evaluate the homogenized acoustic coefficients which are employed in the layer
problem \eq{eq-M1}--\eq{eq-M2}. Its solution, thus, the acoustic pressure in
the waveguide is depicted in Fig.~\ref{fig:num-macro-geom}. Because the flow
field is not completely uniform, local flow problems must be solved either in
all quadrature points of the discretized domain $\Gamma_0$ or in all FE
elements when assuming them constant in elements. By the consequence, the
homogenized coefficients need to be evaluated at the same positions. The
acoustic pressure distribution, obtained for an incident wave with amplitude
$\tilde{p} = 300$\,Pa imposed on $\Gamma_{in}$ and for an anechoic condition
applied on $\Gamma_{out}$, is shown in Fig.~\ref{fig:num-pressure}.

\begin{figure}[ht]
    \centering
    \includegraphics[width=0.98\linewidth]{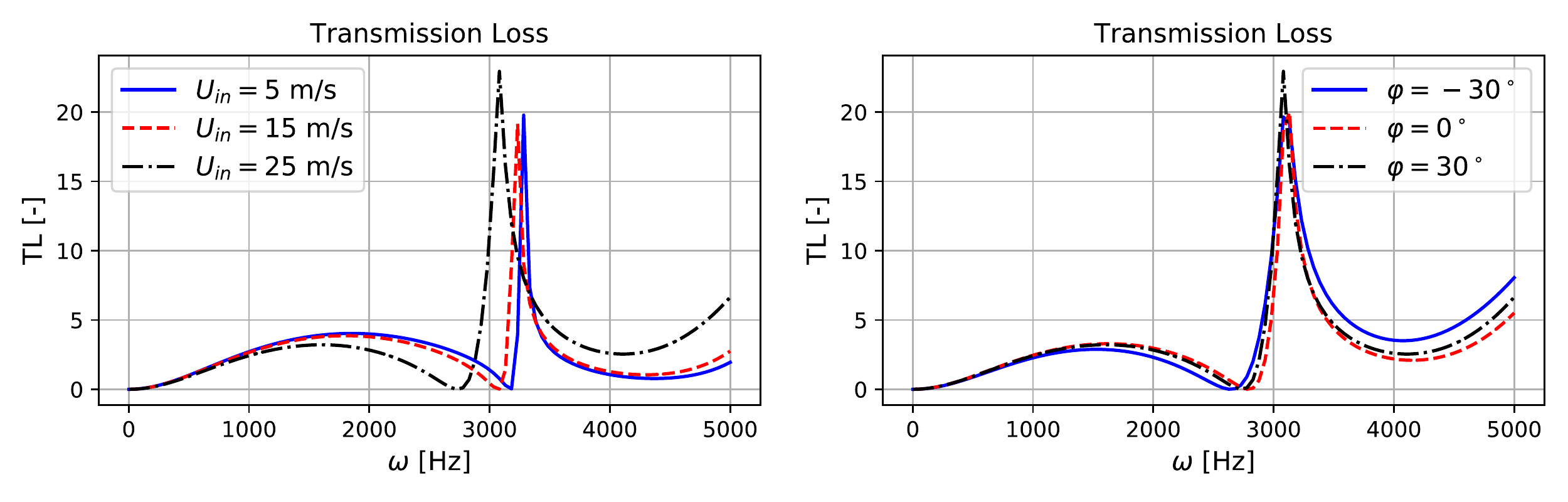}
    \caption{
        Left -- transmission loss (TL) curves calculated for $\varphi=30^\circ$ and
        given input and output velocities: $U_{in} = -U_{out} = 5, 15, 25$ m/s;
        right -- TL curves for $U_{in} = -U_{out} = 25$ m/s and $\varphi=-30, 0, 30^\circ$.}
        
    \label{fig:num-tl}
\end{figure}

The dependence of the global acoustic field on the flow velocity is shown in
Fig.~\ref{fig:num-tl}~left where the transmission loss ${\rm TL}(\omega) =
10\log_{10} \int_{\Gamma_{out}} \vert p \vert^2 / \int_{\Gamma_{in}} \vert p
\vert^2$ is compared for different input and output velocities prescribed on
boundaries $\Gamma_{in}$ and $\Gamma_{out}$: $U_{in} = -U_{out} = 5, 15, 25$
m/s. The global response depends also on the geometrical arrangement of the
perforations which is characterized by the hole slope $\varphi$ in our example.
The effect of various $\varphi$ on the transmission loss is illustrated in
Fig.~\ref{fig:num-tl}~right.

\section{Conclusion}\label{sec-concl}
The derived model of the homogenized acoustic transmission layer serves coupling conditions for the acoustic pressure fields in domains separated by the perforated plate. As the new contribution, we extended the model reported in  \cite{rohan-lukes-waves07} to describe the acoustic waves in flowing fluid. We explored the flow influence on the homogenized coefficients of the homogenized interface model, namely the acoustic impedance. It has been demonstrated, how the flow affects the transmission loss coefficient characterizing the acoustic response of waveguides equipped by the perforated plates.
 Although the background flow increases the computational expenses of the homogenized interface model, as compared to the no-flow case, it enables for a remarkable reduction of the computational effort associated with the numerical simulations of the acoustic field. For small variations of the flow on the interface, the sensitivity analysis of the homogenized coefficients will enable to decrease further the number of the local problems needed to provide a sufficiently accurate numerical approximation of the homogenized flow-dependent interface. Although, in the numerical examples, the potential flow model was employed, the further research is aimed to employ a viscous flow model.

\paragraph{Acknowledgement}
The research has been supported by the grant project GA~2116406S of the Czech Science Foundation.

\appendix

\section{Homogenized coefficients}\label{sec-appx1}

Here we specify the homogenized coefficients (HC) identified in \eq{eq-f6}, when deriving the macroscopic equation \eq{eq-M1} with a nonvanishing  test function $q^0$, and in \eq{eq-cc3}. 
The expressions for the HC are obtained by collecting all terms coupling specific unknown and test functions, as listed below. The symmetry relationships \eq{eq-M3} are proved in \eq{eq-F1}, \eq{eq-M10} and \eq{eq-F7b} using the local autonomous problems \eq{eq-lp3}-\eq{eq-lp5}.

Terms coupling $\gradplx p^0$ and $\gradplx q^0$, where the symmetric expression is due to \eq{eq-lp3},
\begin{equation}\label{eq-F1}
  \begin{split}
     A_{\alpha\beta} & = \ipYs{\nabla_y (y_\alpha + \pi^\alpha)}{\nabla_y y_\beta}
     - \frac{\tau}{c^2} \intY_{Y^*} \wb\cdot \gradply (y_\alpha + \pi^\alpha) \wb\cdot\gradply y_\beta \\
     & =  \ipYs{\nabla_y (y_\alpha + \pi^\alpha)}{\nabla_y(y_\beta + \pi^\beta)}
- \frac{\tau}{c^2} \ipYs{\pd_w^y  (y_\alpha + \pi^\alpha)}{\pd_w^y (y_\beta + \pi^\beta)}
     \;.
\end{split}
\end{equation}

Terms coupling $\gradplx p^0$ and $q^0$,
\begin{equation}\label{eq-F2}
    \imu\om  q^0 \bar\Wb\cdot\gradplx p^0 =  \imu\om\theta  q^0 \intY_{Y^*} \wb\cdot\nabla_y (y_\beta + \pi^\beta) \pd_\beta^x p^0\;, \quad
    \bar W_\beta =  \theta\intY_{Y^*} \wb\cdot\nabla_y (y_\beta + \pi^\beta)\;.
\end{equation}

Terms coupling $p^0$ and $\gradplx q^0$,
\begin{equation}\label{eq-F3}
  \begin{split}
      \imu\om p^0 \bar\Qb^w \gradplx q^0 =
  \imu\om   p^0 \left(c_f^2\intY_{Y^*} \gradply \pi^P - \tau \intY_{Y^*} (\wb\cdot\nabla_y\pi^P)\bar\wb - \theta\intY_{Y^*} \ol\wb     \right) \gradplx q^0
  \;,\\
  \bar \Qb^w = c_f^2\intY_{Y^*} \gradply \pi^P
 - \theta\intY_{Y^*} \bar\wb  - \tau\intY_{Y^*} (\wb\cdot\nabla_y\pi^P)\ol\wb \;,
\end{split}
\end{equation}
where $\ol\wb = (w_\alpha)$, $\alpha = 1,2$ is the in-plane restriction of the advection velocity.


Terms coupling $p^0$ and $q^0$,
\begin{equation}\label{eq-F4}
 \imu\om q^0 M_w p^0 = \imu\om\theta q^0 \intY_{Y^*}\wb\cdot\nabla_y \pi^P p^0\;,\quad
  M_w = \theta  \intY_{Y^*}\wb\cdot\nabla_y \pi^P\;.
\end{equation}

Terms coupling $g^0$ and $\gradplx q^0$,
\begin{equation}\label{eq-F5}
  \begin{split}
    \imu\om c_f^2\gradplx q^0 \Bb g^0 =  \imu\om   \gradplx q^0\cdot\left(c_f^2\intY_{Y^*}\nabla_y \xi\cdot\nabla_y y  - \tau \intY_{Y^*}\ol\wb\otimes\wb\cdot \nabla_y \xi\right) g^0\;,\\
    \Bb =  \intY_{Y^*}\nabla_y \xi\cdot\nabla_y y -\frac{\tau}{c_f^2}\intY_{Y^*}\ol\wb\otimes\wb\cdot \nabla_y \xi 
    \;.
\end{split}
\end{equation}

Terms coupling $g^0$ and $q^0$,
\begin{equation}\label{eq-F6}
    -\om^2 \theta q^0 T_w g^0 =  -\om^2 \theta q^0 \intY_{Y^*}\wb\cdot \nabla_y \xi g^0 \;,\quad
     T_w = \intY_{Y^*}\wb\cdot \nabla_y \xi \;.
\end{equation}

The coupling equation \eq{eq-M2} involves the following homogenized coefficients,
\begin{equation}\label{eq-M10}
\begin{split}
F & = - \intY_{I_y^+} \xi + \intY_{I_y^-} \xi\;, 
\\
  {B'}_\alpha & = \intY_{I_y^+} \pi^\alpha - \intY_{I_y^-}  \pi^\alpha =
 \ipYf{\AA\xi}{\pi^\alpha} = -
  \intY_{Y_\vkappa^*} \pd_\alpha^y \xi +  \frac{\tau}{c_f^2}\intY_{Y^*} w_\alpha \pd_\alpha^y \xi
  = B_\alpha\;,\quad \alpha = 1,2\;,\\
T_w' & =\intY_{I_y^+} \pi^P - \intY_{I_y^-}  \pi^P = - \ipYf{\AA\pi^P}{\xi} = - \frac{\theta}{c_f^2}\intY_{Y^*}\wb\cdot\nabla_y\xi = - \frac{\theta}{c_f^2} T_w\;,
\end{split}
\end{equation}
where \eq{eq-lp4}-\eq{eq-lp5} was employed 

Finally, to show 
$\bar\Wb = -\bar \Wb'$, we employ the following identity obtained due to \eq{eq-lp3} and \eq{eq-lp5},
\begin{equation}\label{eq-F7a}
  \begin{split}
    c_f^2\intY_{Y^*} \gradply \pi^P - \tau \intY_{Y^*} (\wb\cdot\nabla_y\pi^P)\bar\wb =
     c_f^2\ipYs{\gradply \pi^P}{\gradply y_\beta} - \tau \ipYs{\pd_w^y y_\beta}{\pd_w^y \pi^P}\\ = -c_f^2 \ipYf{\AA\pi^\beta}{\pi^P} = -\theta c_f^2\intY_{Y^*}\wb\cdot\nabla_y \pi^\beta\;,
\end{split}
\end{equation}
hence
\begin{equation}\label{eq-F7b}
  \begin{split}
   \bar W_\beta' = c_f^2\intY_{Y^*} \gradply \pi^P
   - \tau\intY_{Y^*} (\wb\cdot\nabla_y\pi^P)\ol w_\beta  - \theta\intY_{Y^*} \bar w_\beta \\
    = \theta\intY_{Y^*}\wb\cdot\nabla_y \pi^\beta- \theta\intY_{Y^*} \bar\wb =
   -\theta\intY_{Y^*} \wb\cdot\nabla_y (y_\beta + \pi^\beta) = - \bar W_\beta\;.
\end{split}
\end{equation}

\bibliographystyle{elsarticle-num}
\bibliography{biblio-acoustics}

\end{document}